\documentclass[12pt,leqno,draft]{article}

\begin{document}

\title{Some remarks about Cauchy integrals and totally real surfaces in
${\bf C}^m$}

\author{Stephen Semmes}

\date{}

\maketitle

	Let us begin by reviewing some geometrically-oriented linear
algebra, about which Reese Harvey once tutored me.  Fix a positive
integer $m$.  The standard Hermitian inner product on ${\bf C}^m$
is defined by
\begin{equation}
	\langle v, w \rangle = \sum_{j=1}^m v_j \, \overline{w_j},
\end{equation}
where $v$, $w$ are elements of ${\bf C}^m$ and $v_j$, $w_j$ denote their
$j$th components, $1 \le j \le m$.  This expression is complex-linear
in $v$, conjugate-complex-linear in $w$, and satisfies
\begin{equation}
	\langle w, v \rangle = \overline{\langle v, w \rangle}.
\end{equation}
Of course $\langle v, v \rangle$ is the same as $|v|^2$, the square
of the standard Euclidean length of $v$.

	Define $(v, w)$ to be the real part of $\langle v, w \rangle$.
This is a real inner product on ${\bf C}^m$, which is real linear in
both $v$ and $w$, symmetric in $v$ and $w$, and such that $(v, v)$ is
also equal to $|v|^2$.  This is the same as the standard real inner
product on ${\bf C}^m \approx {\bf R}^{2m}$.

	Now define $[v, w]$ to be the imaginary part of $\langle v, w \rangle$.
This is a real linear function in each of $v$ and $w$, and it is
antisymmetric, in the sense that
\begin{equation}
	[w, v] = - [v, w].
\end{equation}
Also, $[v, w]$ is nondegenerate, which means that for each nonzero $v$
in ${\bf C}^m$ there is a $w$ in ${\bf C}^m$ such that $[v, w] \ne 0$.
Indeed, one can take $w = i \, v$.

	Let $L$ be an $m$-dimensional real-linear subspace of ${\bf C}^m$.
We say that $L$ is totally-real if $L$ is transverse to $i \, L$,
where $i \, L = \{ i \, v : v \in L\}$.  Transversality here can be
phrased either in terms of $L \cap i \, L = \{0\}$, or in terms
of $L + i \, L = {\bf C}^m$.

	An extreme version of this occurs when $i \, L$ is the
orthogonal complement of $L$.  Because we are assuming that $L$ has
real dimension $m$, this is the same as saying that elements
of $i \, L$ are orthogonal to elements of $L$.  This is equivalent
to saying that $[v, w] = 0$ for all $v$, $w$ in $L$.  Such a real
$m$-dimensional plane is said to be Lagrangian.

	As a basic example, ${\bf R}^m$ is a Lagrangian subspace
of ${\bf C}^m$.  In fact, the Lagrangian subspaces of ${\bf C}^m$
can be characterized as images of ${\bf R}^m$ under unitary linear
transformations on ${\bf C}^m$.  The images of ${\bf R}^m$ under
special unitary linear transformations, which is to say unitary
transformations with complex determinant equal to $1$, are
called special Lagrangian subspaces of ${\bf C}^m$.

	Now suppose that $M$ is some kind of submanifold or surface in
${\bf C}^m$ with real dimension $m$.  We assume at least that $M$ is a
closed subset of ${\bf C}^m$ which is equipped with a nonnegative
Borel measure $\mu$, in such a way that $M$ is equal to the support of
$\mu$, and the $\mu$-measure of bounded sets are finite.  One might
also ask that $\mu$ behave well in the sense of a doubling condition
on $M$, or even Ahlfors-regularity of dimension $m$.  One may wish to
assume that $M$ is reasonably smooth, and anyway we would ask that $M$
is at least rectifiable, so that $\mu$ can be written as the
restriction of $m$-dimensional Hausdorff measure to $M$ times a
density function, and $M$ has $m$-dimensional approximate tangent
spaces at almost all points.

	Let us focus on the case where $M$ is totally real, so that
its approximate tangent planes are totally real, at least almost
everywhere.  In fact one can consider quantitative versions of this.
Namely, if 
\begin{equation}
	d\nu_m = dz_1 \wedge dz_2 \wedge \cdots \wedge dz_m
\end{equation}
is the standard complex volume form on ${\bf C}^m$, then a linear
subspace $L$ of ${\bf C}^m$ of real dimension $m$ is totally real
if and only if the restriction of $d\nu_m$ to $L$ is nonzero.
In any event, the absolute value of the restriction of $d\nu_m$
to $L$ is equal to a nonnegative real number times the standard
positive element of $m$-dimensional volume on $L$, and positive
lower bounds on that real number correspond to quantitative 
measurements of the extent to which $L$ is totally real.
In the extreme case when $L$ is Lagrangian, this real number
is equal to $1$.  For the surface $M$, one can consider lower
bounds on this real coefficient at each point, or at least 
almost everywhere.

	From now on let us assume that $M$ is oriented, so that the
approximate tangent planes to $M$ are oriented.  This means that
reasonably-nice complex-valued functions on $M$ can be integrated
against the restriction of $d\nu_m$ to $M$.  One can then define
pseudo-accretivity and para-accretivity conditions for the restriction
of $d\nu_m$ to $M$ as in \cite{D-J-S}, which basically mean that
classes of averages of the restriction of $d\nu_m$ to $M$ have nice
lower bounds for their absolute values compared to the corresponding
averages of the absolute value of the restriction to $d\nu_m$ to $m$.
This takes into account the oscillations of the restriction of $d\nu_m$
to $M$.

	Note that if $M$ is a smooth submanifold of ${\bf C}^m$ of
real dimension $m$, then $M$ is said to be Lagrangian if its tangent
spaces are Lagrangian $m$-planes at each point.  This turns out to be
equivalent to saying that $M$ can be represented locally at each point
as the graph of the gradient of a real-valued smooth function on ${\bf
R}^m$ in an appropriate sense, as in \cite{weinstein}.  If the tangent
planes of $M$ are special Lagrangian, then $M$ is said to be a special
Lagrangian submanifold.  See \cite{reese, reese-blaine} in connection
with these.

	It seems to me that there is a fair amount of room here for
various interesting things to come up, basically concerning the
geometry of $M$ and aspects of several complex variables on ${\bf
C}^m$ around $M$.  When $m = 1$, this would include the Cauchy
integral operator applied to functions on a curve and holomorphic
functions on the complement of the curve.  In general this can include
questions about functional calculi, as in \cite{C-M1, C-M2, D-J-S},
and $\overline{\partial}$ problems with data of type $(0,m)$, as well
as relations between the two.


\begin{thebibliography}{99}


\bibitem {A-H-L2-M-T} P.~Auscher, S.~Hofmann, M.~Lacey, J.~Lewis,
A.~McIntosh, and P.~Tchamitchian, {\it The solution of Kato's
conjectures}, Comptes Rendus des Sc\'eances de l'Acad\'emie des
Sciences de Paris S\'er.\ I {\bf 322} (2001), 601--606.

\bibitem {A-H-L-M-T} P.~Auscher, S.~Hofmann, M.~Lacey, A.~McIntosh,
and P.~Tchamitchian, {\it The solution of the Kato square root problem
for second order elliptic operators on ${\bf R}^n$}, Annals of
Mathematics (2) {\bf 156} (2002), 633--654.

\bibitem {A-H-L-T} P.~Auscher, S.~Hofmann, J.~Lewis, and
P.~Tchamitchian, {\it Extrapolation of Carleson measures and the
analyticity of Kato's square-root operators}, Acta Mathematica {\bf
187} (2001), 161--190.

\bibitem {A-H-M-T} P.~Auscher, S.~Hofmann, A.~McIntosh, and
P.~Tchamitchian, {\it The Kato square root problem for higher order
elliptic operators and systems on ${\bf R}^n$, Dedicated to the memory
of Tosio Kato}, Journal of Evolution Equations {\bf 1} (2001),
361--385.

\bibitem {A-T} P.~Auscher and P.~Tchamitchian, {\it Square Root
Problem for Divergence Operators and Related Topics}, Ast\'erisque
{\bf 249}, 1998.

\bibitem {bell} S.~Bell, {\it The Cauchy Transform, Potential Theory,
and Conformal Mapping}, CRC Press, 1992.

\bibitem {C-M-M} R.~Coifman, A.~McIntosh, and Y.~Meyer, {\it
L'Int\'egrale de Cauchy d\'efinit un op\'erateur born\'e sur $L^2$
pour les courbes lipschitziennes}, Annals of Mathematics (2) {\bf 116}
(1982), 361--387.

\bibitem {C-M1} R.~Coifman and Y.~Meyer, {\it Fourier analysis of
multilinear convolutions, Calder\'on's theorem, and analysis on
Lipschitz curves}, in {\it Euclidean Harmonic Analysis}, Proceedings
of Seminars held at the University of Maryland, 1979, 104--122,
Lecture Notes in Mathematics {\bf 779}, Springer-Verlag, 1980.

\bibitem {C-M2} R.~Coifman and Y.~Meyer, {\it Non-linear harmonic
analysis, operator theory, and PDE}, in {\it Beijing Lectures in
Harmonic Analysis}, Annals of Mathematics Studies {\bf 112}, 3--45,
Princeton University Press, 1986.	

\bibitem {C-S1} R.~Coifman and S.~Semmes, {\it Real-analytic
operator-valued functions defined in BMO}, in {\it Analysis and
Partial Differential Equations}, 85--100, Marcel Dekker, 1990.

\bibitem {C-S2} R.~Coifman and S.~Semmes, {\it $L^2$ estimates in
nonlinear Fourier analysis}, in {\it Harmonic Analysis (Sendai,
1990)}, Proceedings of the ICM-90 Satellite Conference, 79--95,
Springer-Verlag, 1991.

\bibitem {C-D-M-Y} M.~Cowling, I.~Doust, A.~McIntosh, and A.~Yagi,
{\it Banach space operators with a bounded $H^\infty$ functional
calculus}, Journal of the Australian Mathematical Society Ser.\ A {\bf
60} (1996), 51--89.

\bibitem {D-J-S} G.~David, J.-L.~Journ\'e, and S.~Semmes, {\it
Op\'erateurs de Calder\'on--Zygmund, fonctions para-accr\'etives et
interpolation}, Revista Matem\'atica Iberoamericana {\bf 1} (4)
(1985), 1--56.

\bibitem {peter} P.~Duren, {\it Theory of $H^p$ Spaces}, Academic
Press, 1970.

\bibitem {john} J.~Garnett, {\it Bounded Analytic Functions},
Academic Press, 1981.

\bibitem {gleason1} A.~Gleason, {\it The abstract theorem of
Cauchy--Weil}, Pacific Journal of Mathematics {\bf 12} (1962),
511--525.

\bibitem {gleason2} A.~Gleason, {\it The Cauchy--Weil theorem},
Journal of Mathematics and Mechanics {\bf 12} (1963), 429--444.

\bibitem {reese} R.~Harvey, {\it Spinors and Calibrations},
Academic Press, 1990.

\bibitem {reese-blaine} R.~Harvey and B.~Lawson, {\it Calibrated
geometries}, Acta Mathematica {\bf 148} (1982), 47--157.

\bibitem {H-L} G.~Henkin and J.~Leiterer, {\it Andreotti--Grauert
Theory by Integral Formulas}, Akademie-Verlag, 1988.

\bibitem {hormander} L.~H\"ormander, {\it An Introduction to Complex
Analysis in Several Variables}, North-Holland, 1973.

\bibitem {J-M-PW} B.~Jefferies, A.~McIntosh, and J.~Picton-Warlow,
{\it The monogenic functional calculus}, Studia Mathematica {\bf 136}
(1999), 99--119.

\bibitem {jean-lin} J.-L.~Journ\'e, {\it Calder\'on--Zygmund Operators,
Pseudo-Differential Operators, and the Cauchy Integral of Calder\'on},
Lecture Notes in Mathematics {\bf 994}, Springer-Verlag, 1983.

\bibitem {K-M} C.~Kenig and Y.~Meyer, {\it Kato's square roots of
accretive operators and Cauchy kernels on Lipschitz curves are the
same}, in {\it Recent Progress in Fourier Analysis (El Escorial,
1983)}, 123--143, North-Holland, 1985.

\bibitem {krantz} S.~Krantz, {\it Function Theory of Several Complex
Variables}, second edition, AMS Chelsea Publishing, 2001.

\bibitem {alan} A.~McIntosh, {\it Operators which have an $H^\infty$
functional calculus}, in {\it Miniconference on Operator Theory and
Partial Differential Equations (North Ryde, 1986)}, 210--231,
Proceedings of the Centre for Mathematical Analysis {\bf 14},
Australian National University, 1986.

\bibitem {M-M} A.~McIntosh and Y.~Meyer, {\it Alg\`ebres
d'op\'erateurs d\'efinis par des int\'egrales singuli\`eres}, Comptes
Rendus des Sc\'eances de l'Acad\'emie des Sciences de Paris S\'er.\ I
Math.\ {\bf 301} (1985), 395--397.

\bibitem {M-P} A.~McIntosh and A.~Pryde, {\it A functional calculus
for several commuting operators}, Indiana University Mathematics
Journal {\bf 36} (1987), 421--439.

\bibitem {M-Y} A.~McIntosh and A.~Yagi, {\it Operators of type
$\omega$ without a bounded $H^\infty$ functional calculus}, in {\it
Miniconference on Operators in Analysis (Sydney, 1989)}, 159--172,
Proceedings of the Centre for Mathematical Analysis {\bf 24},
Australian National University, 1990.

\bibitem {pertti} P.~Mattila, {\it Geometry of Sets and Measures in
Euclidean Spaces: Fractals and Rectifiability}, Cambridge University
Press, 1995.

\bibitem {ss1} S.~Semmes, {\it A criterion for the boundedness of singular
integrals on hypersurfaces}, Transactions of the American Mathematical
Society {\bf 311} (1989), 501--513.

\bibitem {ss3} S.~Semmes, {\it Chord-arc surfaces with small constant I},
Advances in Mathematics {\bf 85} (1991), 198--223.

\bibitem {joan} J.~Verdera, {\it $L^2$ boundedness of the Cauchy
integral and Menger curvature}, in {\it Harmonic Analysis and Boundary
Value Problems (Fayetteville AR, 2000)}, 139--158, Contemporary
Mathematics {\bf 277}, American Mathematical Society, 2001.

\bibitem {weinstein} A.~Weinstein, {\it Lectures on Symplectic
Manifolds}, Conference Board of the Mathematical Sciences Regional
Conference Series in Mathematics {\bf 29}, American Mathematical
Society, 1977.



\end{thebibliography}
\end{document}